 \documentclass[11p,times]{article}

 \usepackage{graphicx}

\usepackage{amssymb}

\begin{document}




\title{Maximum entropy estimation of probability distributions with Gaussian conditions}
\maketitle

\begin{center} \author{Mihail - Ioan Pop \\ Department of Electrical Engineering and Applied Physics, Transilvania University of Bra\c{s}ov, Romania \\ e-mail: mihailp@unitbv.ro} \end{center}

\begin{abstract}
We describe a method to computationally estimate the probability density function of a univariate random variable by applying the maximum entropy principle with some local conditions given by Gaussian functions. The estimation errors and optimal values of parameters are determined. Experimental results are presented. The method estimates the distribution well if a large enough selection is used, typically at least 1 000 values. Compared to the classical approach of entropy maximisation, local conditions allow improving estimation locally. The method is well suited for a heuristic optimisation approach.
\end{abstract}







\section{Description of method}

Consider a continuous random variable $X \in \mathbb{R}$ with probability density function (pdf) $\rho (x)$ and a selection of $N$ values $ \left( X_i \right) $, $i = 1,2,...,N$ of $X$. We assume $\rho (x)$ to be of class $C^2$ everywhere. The purpose of the described method is to computationally estimate $\rho (x)$ using this selection of $X$. For this, $\rho (x)$ is restricted to the interval $[ X_{\min} ; X_{\max} ]$, where $X_{\min} = \min \{ X_i : i = 1,2,...,N \}$, $X_{\max} = \max \{ X_i : i = 1,2,...,N \}$, and label $\Delta X = X_{\max} - X_{\min} $. Next the pdf is discretised on $N_p$ equidistant points $x_j \in [ X_{\min} ; X_{\max} ]$, $j = 1,2,...,N_p$ of the form $x_j = X_{\min} + \Delta X \cdot (j-1) / ( N_p - 1 ) $, generating a probability distribution $\left( p_j \right)$, $j = 1,2,...,N_p$. The values of $p_j$ are computed from the selection $ \left( X_i \right)$ as presented below. Next, the pdf $\rho (x)$ is estimated on each $x_j$ as $\rho \left( x_j \right) \approx p_j \Delta X \cdot (j-1) / ( N_p - 1 )$. 

In order to apply the maximum entropy principle \cite{Jaynes}, some conditions are imposed. These are constructed with a function $f : \mathbb{R} \rightarrow \mathbb{R}$, $f(x) = \exp \left( - x^2 / 2 \right)$, $x \in \mathbb{R}$. These functions are inspired by the Kernel Density Estimation method \cite{Parzen,Rosenblatt}.

We center the function $f (x)$ on different points in the domain of values of $X$. For this, take $N_c$ equidistant values $c_k \in [ X_{\min} ; X_{\max} ]$  of the form $c_k = X_{\min} + \Delta X \cdot (k-1) / ( N_c - 1 ) $ and parameters $\sigma _k > 0$, $k = 1,2,...,N_c$. For each $c_k$ a function $f_k$ is built as $f_k (x) = f \left( \frac{ x - c_k }{ \sigma _k } \right)$. Next, two averages of $f_k$ are computed: an empirical average given by selection values: 
\begin{equation}
F_k ^{emp} = \frac{1}{N} \sum_{i=1}^{N} f_k \left( X_i \right)
\end{equation}
and a simulated average given by estimated probabilities $p_j$: 
\begin{equation}
F_k ^{sim} = \frac{1}{N_p} \sum_{j=1}^{N_p} p_j f_k \left( x_j \right) .
\end{equation}
The probability distribution $( p_j )$ is determined by maximising the Shannon entropy of $( p_j )$: $H = - \sum _{j=1}^{N_p} p_j \ln p_j$ with conditions $F_k ^{sim} = F_k ^{emp}$, $k = 1,2,...,N_c$ and $\sum _{j=1}^{N_p} p_j = 1$. For ease of computation, $\left( p_j \right)$ is obtained by minimising a cost function of the form: 
\begin{equation}
E = \sum_{k=1}^{N_c} \left( F_k ^{sim} - F_k ^{emp} \right)^2 - k_H \cdot H, 
\end{equation}
where $k_H > 0$. This relaxes the above conditions. The unit sum of $\left( p_j \right)$ is imposed at the end by dividing each obtained $p_j$ to $\sum _{j=1}^{N_p} p_j$.

The value of $k_H$ regulates the smoothness of the estimated pdf. This is necessary because the empirical averages $F_k ^{emp}$ contain statistical noise, which is transmitted to the estimated pdf. By maximising the entropy the noise is reduced. Nevertheless, too great an importance given to entropy maximisation may lead to smoothing real features of the pdf. The parameters $\sigma _k$ control the locality and the smoothness of the estimated pdf, \textit{i.e.} a smaller value of $\sigma _k$ define conditions on a smaller neighbourhood of $c_k$, allowing a better local approximation of $\rho \left( c_k \right)$. In practice, the minimisation was done on a computer with the Simulated Annealing (SA) algorithm \cite{Kirk}. Because of the finite number of steps of this algorithm, some more noise is introduced in the final pdf estimation. This is eliminated by applying a moving average on the final result.

\section{Errors of estimation} 

We determine the errors of the pdf estimation and calculate optimal values of $\sigma _k$. We work in a small perturbations setting. Experimentally, the algorithm was observed to reproduce the theoretical pdf well for $N \geq 1\,000$ selection values, which supports this assumption. 

The algorithm converges to a perturbed pdf $\rho _p (x)$ with respect to the theoretical pdf $\rho (x)$, where $\rho, \rho _p : D \rightarrow \mathbb{R}$, $D \subset \mathbb{R}$. The perturbation of the pdf is $\delta \rho (x) = \rho _p (x) - \rho (x)$. Since both $\rho$ and $\rho _p$ have unit integral, it follows that $\int _D \delta \rho (x) dx = 0$.

Consider the case where there is only one condition given by a function $f_1 (x)$ centered at $c_1 = c \in D$ with parameter $\sigma _1 = \sigma$. We are interested in the behaviour of the estimated pdf in a neighbourhood of $c$. Replace $ F_k ^{emp}$ by the theoretical average $F = \int f_1(x) \rho (x) dx$ and $ F_k ^{sim}$ by the perturbed average $F_{p} = \int _D f_1(x) \rho _p (x) dx$. The perturbation of the average is $\delta F = F_p - F = \int _D f_1(x) \delta \rho (x) dx$. Then the cost function is perturbed to $\delta E = E_p - E = - k_H \delta H + \left( \delta F \right)^2$. 

Here we take the continuous Shannon entropy $H = - \int _D \rho (x) \ln \rho (x) dx$. We consider that $| \delta \rho (x) | \ll \rho (x)$, $\forall x \in \mathbb{R}$. Then we have from the Taylor series: 
\begin{equation}
\rho _p (x) \ln \rho _p (x) = \rho (x) \ln \rho (x) + \delta \rho (x) \ln \rho (x) + \delta \rho (x) + \frac{1}{2} \frac{ \left( \delta \rho (x) \right) ^2}{\rho (x)} + O \left( \delta \rho \right) ^3
\end{equation}
as $\delta \rho (x) \rightarrow 0$. We compute the entropy perturbation to the second degree with respect to $\delta \rho (x)$.
The perturbation of the cost function becomes: 
\begin{equation}
\delta E = k_H \left( \int _D \delta \rho (x) \ln \rho (x) dx + \frac{1}{2} \int _D \frac{ \left( \delta \rho (x) \right) ^2}{\rho (x)} dx \right) + \left( \int _D f(x) \delta \rho (x) dx \right) ^2 . 
\end{equation}
The condition $E _p = \min$ reduces to $\delta E = \min$, since $E$ is fixed by $\rho$. Then, by the method of Lagrange multipliers, $\delta \rho (x)$ can be determined by minimising a functional $\delta \Lambda = \delta E + \lambda \int _D \delta \rho (x) dx$, $\lambda \in \mathbb{R}$. For a small variation $\delta ' \delta \rho (x)$ of $\delta \rho (x)$ there corresponds a small variation $\delta ' \delta \Lambda$ of $\delta \Lambda$ and the condition of minimum reduces to $\delta ' \delta \Lambda = 0$, \textit{i.e.}: 
\begin{equation}
\int _D \left( k_H \frac{\delta \rho (x)}{\rho (x)} + k_H \ln \rho (x) + 2 f_1(x) \delta F + \lambda \right) \delta ' \delta \rho (x) dx = 0 . 
\end{equation}
It holds for any $\delta ' \delta \rho (x)$ if and only if the paranthesis under the integral is zero. The multiplier $\lambda$ is obtained from $\int _D \delta \rho (x) dx = 0$. The perturbation of the pdf becomes:
\begin{equation} \label{Sol}
\delta \rho (x) = - \rho (x) \left[ \frac{2}{k_H} \left( f_1(x) - F \right) \delta F + \ln \rho (x) + H \right] .
\end{equation}

We label the variation of $\rho$ around $c$ $\Delta \rho (x) = \rho (x) - \rho (c)$. Now consider the interval $D$ is chosen such that $| \Delta \rho (x) | \ll \rho (x)$, $\forall x \in D$. This assumption is not very restrictive, since $D$ can be taken as a subset of the interval of values of $X$ in order to study $\rho (x)$ locally. We label $\Delta P = \int _D \Delta \rho (x) dx$. Then $H$ can be further expressed by Taylor series development of $\rho \ln \rho$ as $H = - \ln \rho (c) - \Delta P + O \left( \Delta \rho \right) ^2$ and $\ln \rho (x) = \ln \rho (c) + \Delta \rho (x) / \rho (c) + O \left( \Delta \rho \right) ^2$, such that (\ref{Sol}) becomes: 
\begin{equation}
\delta \rho (x) = - \rho (x) \left[ \frac{2}{k_H} \left( f_1(x) - F \right) \delta F + \frac{\Delta \rho (x)}{\rho (c)} - \Delta P \right] + O \left( \Delta \rho \right) ^2 .
\end{equation}
Here $\delta \rho (x)$ is identified with the error of the estimation of pdf, while $\delta F$ represents the error of the average of $f_1$. The free term $\Delta \rho (x) / \rho (c) - \Delta P$ is an error given by the local variation of the pdf around $x = c$. We have $\Delta \rho (c) = 0$ and $f_1 (c) = 1$. 

The total error of conditions $\delta F ^{tot}$ has two sources: the computation error of $F_1^{emp}$, which gives a term $\delta F ^{emp}$ and the pdf estimation process, which gives a term $\delta F ^{est}$. The term $\delta F ^{emp}$ can be estimated by considering that $F_1^{emp}$ is computed with a large enough number $N$ of values from the selection $ \left( X_i \right)$, such that $F_1^{emp}$ is approximately normally distributed. From the law of large numbers $ \delta F_1^{emp} $ is a random variable with approximately normal distribution with average $E \delta F ^{emp} = 0$ and variance $Var \delta F ^{emp} = Var (f_1(X)) / N$, 
where $Var \left( f_1(X) \right)$ is the variance of $f_1$: $Var \left( f_1(X) \right) = F_2 - F^2$, with $F_2 = \int _D f_1 ^2 (x) \rho (x) dx$.

We identify $\delta F ^{est} = \delta F = \int _D f_1 (x) \delta \rho (x) dx$. Label $H_1 = - \int _D f_1 (x) \rho (x) \ln \rho (x) dx$ and $\Delta F = \int _D f_1 (x) \Delta \rho (x) dx$. Then the expression of $\delta \rho (x)$ becomes: 
\begin{equation}
\delta \rho (x) = - \rho (x) \left[ \frac{2}{k_H} \left( f_1(x) - F \right) \left( \delta F ^{emp} + \delta F \right) + \ln \rho (x) + H \right] .
\end{equation}
Putting this expression into the integral form of $\delta F$, we get: 
\begin{equation}
\delta F = \delta F ^{est} = \frac{H_1 - H F - \frac{2}{k_H} Var (f_1(X)) \cdot \delta F ^{emp}}{ 1 + \frac{2}{k_H} Var (f_1(X)) } . 
\end{equation}
If $| \Delta \rho (x) | \ll \rho (x)$, $\forall x \in D$, by the Taylor series of $\rho \ln \rho$ we have $H_1 - H F = F \Delta P - \Delta F + O (\Delta \rho )^2 $ and the error of the pdf estimation becomes: 
\begin{equation} \label{eroare}
\delta \rho (c) = - \frac{2}{k_H} \rho (c) \left( 1 - F \right) \frac{ \delta F ^{emp} + F \Delta P - \Delta F }{ 1 + \frac{2}{k_H} Var (f_1(X)) } + \rho (c) \Delta P . 
\end{equation}
It follows that, for large enough $N$, $\delta \rho (c)$ is an approximately normally distributed random variable with parameters:
\begin{eqnarray}
E \delta \rho (c) &=& - \frac{2}{k_H} \rho (c) \left( 1 - F \right) \frac{ F \Delta P - \Delta F }{ 1 + \frac{2}{k_H} Var (f_1(X)) } + \rho (c) \Delta P , \\
Var \delta \rho (c) &=& \frac{4}{N k_H^2} \rho ^2 (c) \left( 1 - F \right) ^2 \frac{ Var (f_1(X)) }{ \left( 1 + \frac{2}{k_H} Var (f_1(X)) \right) ^2 } .
\end{eqnarray}
In the same way, the total error of conditions $\delta F ^{tot}$ is a normally distributed random variable
\begin{equation}
\delta F ^{tot} = \frac{F \Delta P - \Delta F + \delta F ^{emp}}{ 1 + \frac{2}{k_H} Var (f_1(X)) } 
\end{equation}
with parameters:
\begin{eqnarray}
E \delta F ^{tot} &=& \frac{F \Delta P - \Delta F }{ 1 + \frac{2}{k_H} Var (f_1(X)) } , \\
Var \delta F ^{tot} &=& \frac{1}{N} \frac{ Var (f_1(X)) }{ \left( 1 + \frac{2}{k_H} Var (f_1(X)) \right) ^2 } .
\end{eqnarray}

Take $D = [c-d/2;c+d/2]$ and $\rho (x) \approx \rho (c) + \rho ' (c) (x-c) + \frac{\rho '' (c)}{2} (x-c)^2$, $x \in D$. With the Gaussian function $f_1 (x)$ we have:
\begin{eqnarray}
F &=& \sqrt{2\,\pi }\,\left( \rho (c)\,C_1\,\sigma +\frac{\rho '' (c)}{2}\,C_2\,{\sigma}^{3}\right) , \\
F_2 &=& \sqrt{\pi }\,\left( \rho (c)\,C_3\,\sigma +\frac{\rho '' (c)}{2}\,C_4\,{\sigma}^{3}\right) ,
\end{eqnarray}
where $C_1, C_2, C_3, C_4 \in (0;1]$ are functions of $\sigma$. If $\sigma \ll d$ then $C_i \approx 1$, $i = 1,2,3,4$. 
Also $\Delta P = \rho '' (c)\,\frac{{d}^{3}}{24}$ and $F \Delta P - \Delta F = \frac{\rho '' (c)}{2}\,\left( F \frac{{d}^{3}}{12}-{\sigma}^{3}\,\sqrt{2\,\pi}\,C_2 \right) $.

\section{Pdf error minimisation}

We want to find the value of $\sigma$ that minimises $| \delta \rho (c) |$. Putting the condition $ E \delta \rho (c) = 0$ one obtains a sixth degree equation in $\sigma$, which can be solved numerically only. Nevertheless, this condition can be relaxed for small $k_H$ and $Var (f_1(X))$ to cancellation of the first term in $ E \delta \rho (c)$. The average error of approximation will be then $ E \delta \rho (c) = \rho (c) \Delta P = \rho (c) \rho '' (c) d^3 / 24$. We put the condition $F = 1$, which leads to the third-degree equation, which can be solved exactly \cite{Abramowitz,Press}:
\begin{equation}
\frac{\rho '' (c)}{2} C_2\,{\sigma}^{3} + \rho (c)\,C_1\,\sigma - \frac{1}{\sqrt{2 \pi}} = 0.
\end{equation} 
Its discriminant is $\Delta = - \rho '' (c) \left( \frac{27}{8 \pi} \rho '' (c) C_2^2 + 2 \rho ^3 (c) C_1^3 C_2 \right)$. The parenthesis cancels for $\rho '' (c) = \rho '' _0 (c) = - \frac{16 \pi}{27} \frac{C_1^3}{C_2} \rho ^3 (c)$. If $\rho '' (c) \in \left[\rho '' _0 (c) ; 0 \right]$ then $\Delta \geq 0$ and the equation has 3 real solutions; otherwise, $\Delta < 0$ and the equation has 1 real and 2 complex solutions. If $\rho '' (c) = 0$ there is only one solution: 
\begin{equation}
\sigma _0 = \frac{1}{\sqrt{2 \pi} C_1 \rho (c)} . 
\end{equation}
For $\rho '' (c) \neq 0$ and $\Delta < 0$ label $A = \left[ \frac{1}{\sqrt{2 \pi} C_2 \rho '' (c)} \left( 1 + \sqrt{ 1 - \frac{\rho '' _0 (c)}{\rho '' (c)} } \right) \right] ^{ 1/3 }$. Then the three solutions are:
\begin{eqnarray}
\sigma _1 &=& A - \frac{2 C_1 \rho (c)}{3 A C_2 \rho '' (c)}, \\
\sigma _2 &=& \exp \left( \frac{2 \pi}{3} i \right) A - \exp \left( - \frac{2 \pi}{3} i \right) \frac{2 C_1 \rho (c)}{3 A C_2 \rho '' (c)}, \\
\sigma _3 &=& \exp \left( - \frac{2 \pi}{3} i \right) A - \exp \left( \frac{2 \pi}{3} i \right) \frac{2 C_1 \rho (c)}{3 A C_2 \rho '' (c)}. 
\end{eqnarray}

The first solution is always real. It can be approximated as: $\sigma _1 \approx \sigma _0 - \frac{C_2 \rho '' (c) }{2 {\sqrt{2 \pi}} ^3 C_1^4 \rho ^4 (c) }$ for small $ \rho '' (c) > 0$. The 3 solutions are represented in Figure 1. Only $\sigma > 0$ is meaningful.

\section{Error of conditions minimisation}

Another way to find optimal $\sigma$ is to put the condition $E \delta F ^{tot} = 0$ or equivalently $F \Delta P - \Delta F = 0$. For $D$ centered in $c$ as above, this yields the equation:
\begin{equation}
\sigma \left[C_2\,{\sigma}^{2}\,\left( \rho '' (c)\,{d}^{3}-24\right) +2\,C_1\,\rho (c)\,{d}^{3} \right] = 0 .
\end{equation}
Its only positive solution is:
\begin{equation}
\sigma _4 = \sqrt{\frac{2\,C_1\,\rho (c)\,d^3}{C_2 \left( 24 - \rho '' (c)\,{d}^{3} \right)}} .
\end{equation}
This solution increases with $\rho '' (c)$ and has a vertical asymptote for $\rho '' _1 (c) = 24 / d^3$. 

\section{Experimental verifications and conclusions} 

In practice the conditions were given by Gaussian functions. 
The width $\sigma$ was chosen the same for all $c_k$. Gaussian conditions allowed measuring the error of estimation locally through the relative error of conditions: 
\begin{equation}
\varepsilon \left( c_k \right) = \frac{ \left| F_k ^{sim} - F_k ^{emp} \right|}{ F_k ^{emp} }, k = 1,2,...,N_c . 
\end{equation}
We used $N_c = 101$ conditions and $k_H = 10 ^{-3}$. The pdf was estimated on $N_p = 1\,000$ points and the final pdf was further smoothed with a 10 - point moving average. Tests were carried on the pdf $\rho (x) = \frac{1}{Z} \exp [ - 10^4 \cdot (x-0.1)(x-0.2)(x-0.3)(x-0.5)(x-0.8)(x-0.9) ]$, $x \in [0;1]$, with $Z$ such that $\int _0^1 \rho (x) dx = 1$. Computed values of $\sigma$ vary from 0 up to $\approx \Delta X / 20$ for $\sigma _1$ and up to $\approx \Delta X / 35$ for $\sigma _4$. Results are shown in Figure 1. For $N = 100$ the choice of $\sigma$ is important for the quality of the estimation. Generally, $\sigma _4$ can be used around pdf maxima, while $\sigma _1$ can be used between maxima. For small $| \delta \rho '' (c) |$ $\sigma _0$ may be used. For $N = 1\,000$ the estimation is quite good. The error $\varepsilon \left( c_k \right)$ is small except where $\rho (x) \approx 0$. 

Compared to the classical Maximum Entropy Method using power law conditions, the presented method has some advantages: (i) local conditions allow improving the estimation locally and also measuring its quality with $\varepsilon \left( c_k \right)$; (ii) the discretised pdf is well suited for a heuristic optimisation approach such as the SA, even for high $N_p$, because there are no intrinsic parameters to be determined (the Lagrange multipliers of the classical approach), which the SA has difficulty finding; (iii) the pdf estimate lies between a uniform pdf, obtained for $k_H \rightarrow \infty$ and the experimental values pdf, obtained for $k_H = 0$ and $\sigma \rightarrow 0$; a bad pdf estimate obtained for too small $k_H$ or $\sigma$ can be improved by smoothing; (iv) for large $N$ the method works well for $\sigma$ far from optimal. 
The method was applied to the study of some asteroid parameters \cite{Pop1} and solar cycles \cite{Pop2}.

\begin{figure}
  \centering
  {
	\includegraphics[width=.25\textwidth]{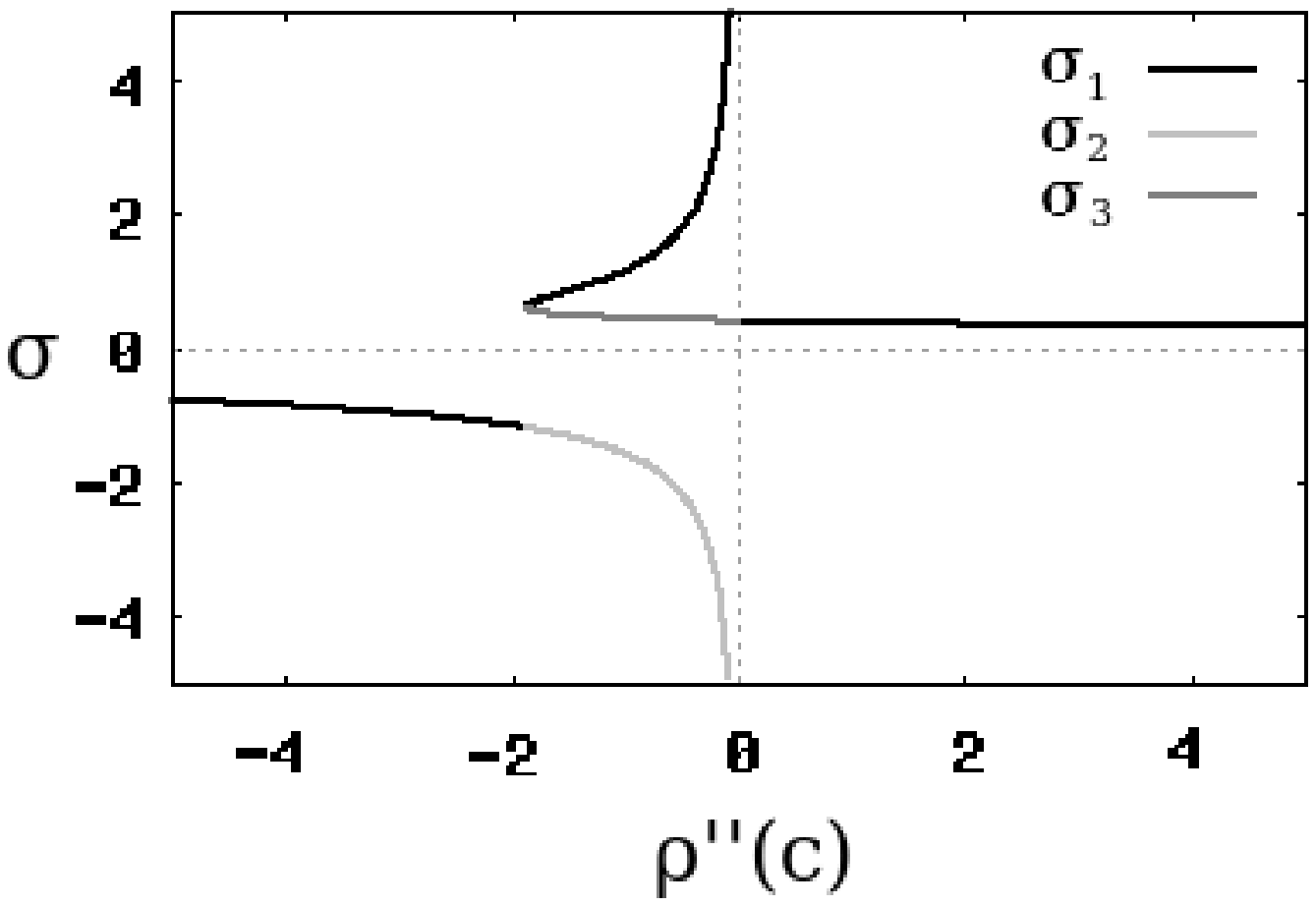}
  	\includegraphics[width=.25\textwidth]{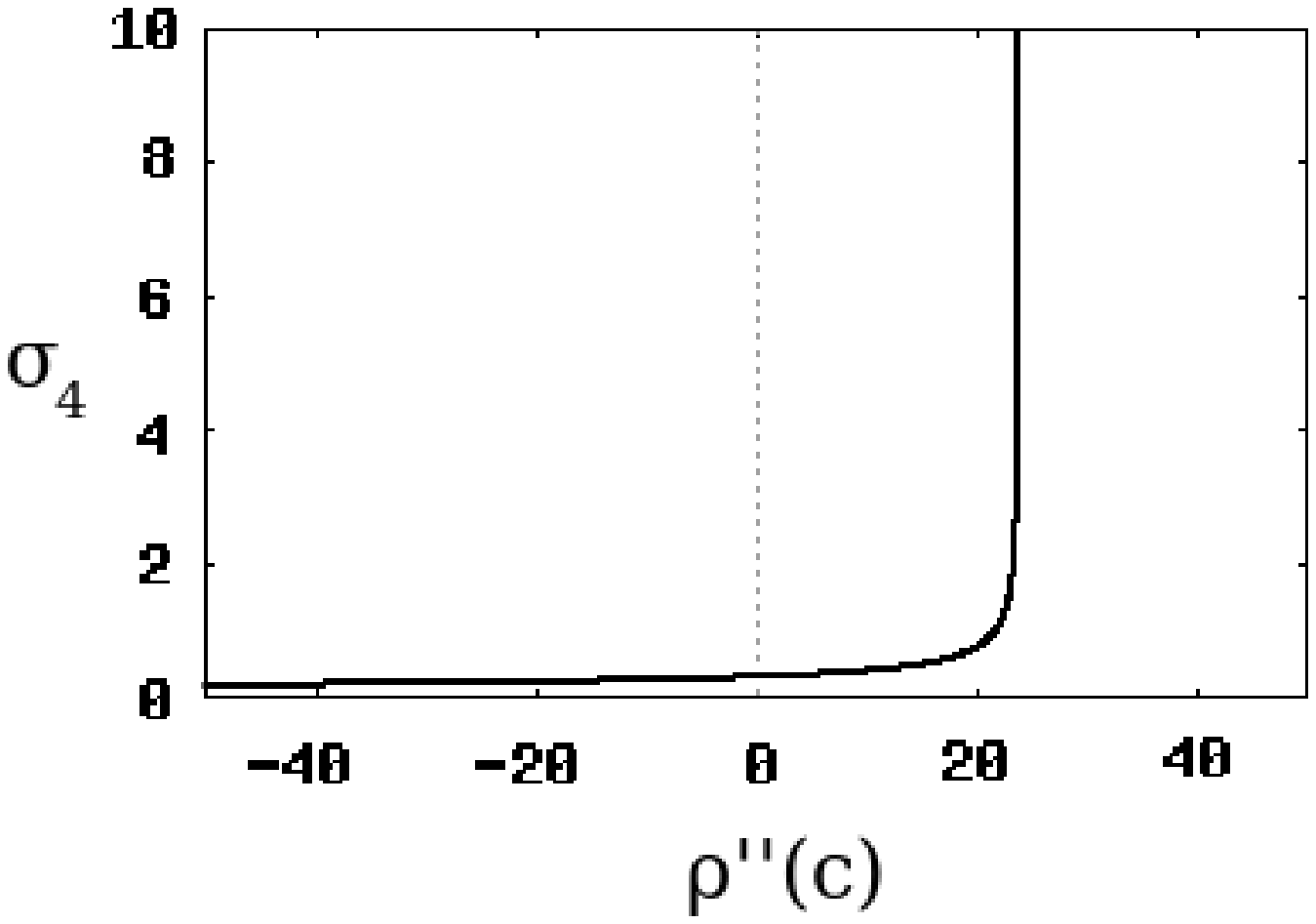}
  }
  \caption{ Variation of optimal $\sigma$ with $ \rho '' (c) $ for $C_1 = C_2 = 1$, $d = 1$, $\rho (c) = 1 $: (left) from pdf error minimisation and (right) from error of conditions minimisation.}
\end{figure}

\begin{figure}
  \centering
  {
	\includegraphics[width=.24\textwidth]{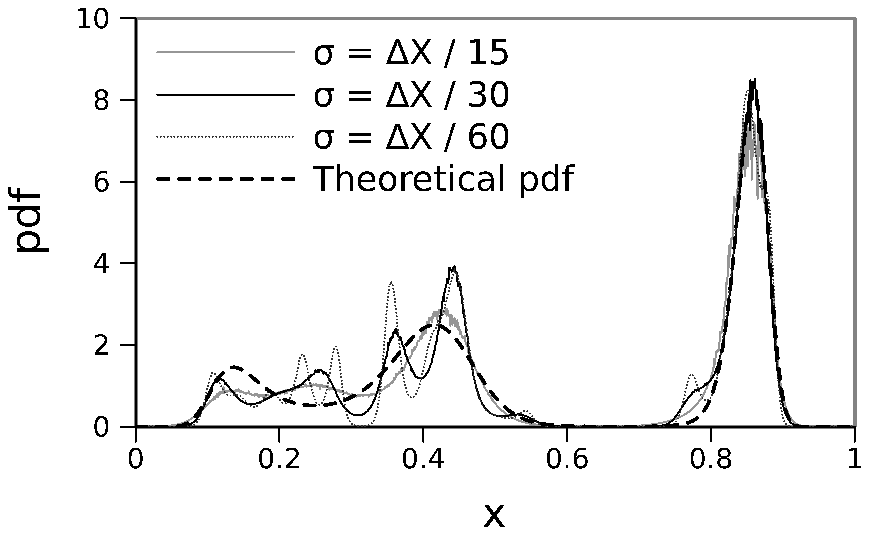}
	\includegraphics[width=.24\textwidth]{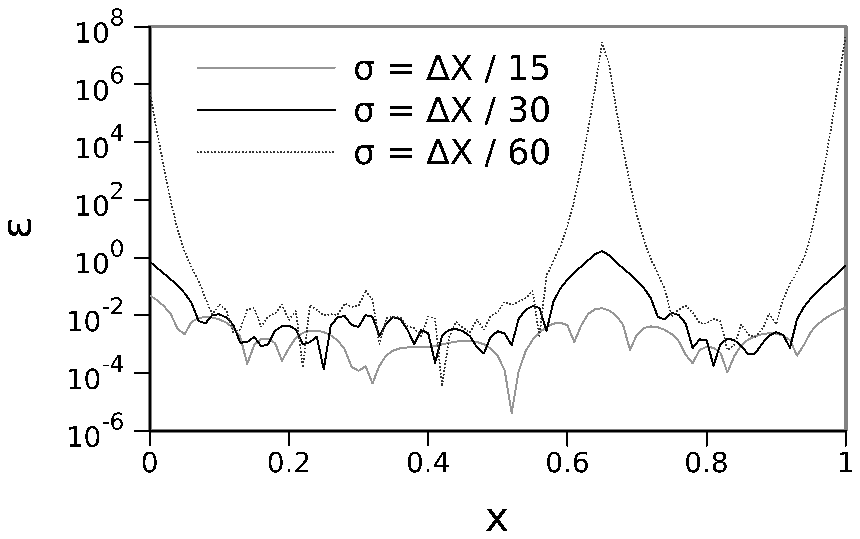}
	\includegraphics[width=.24\textwidth]{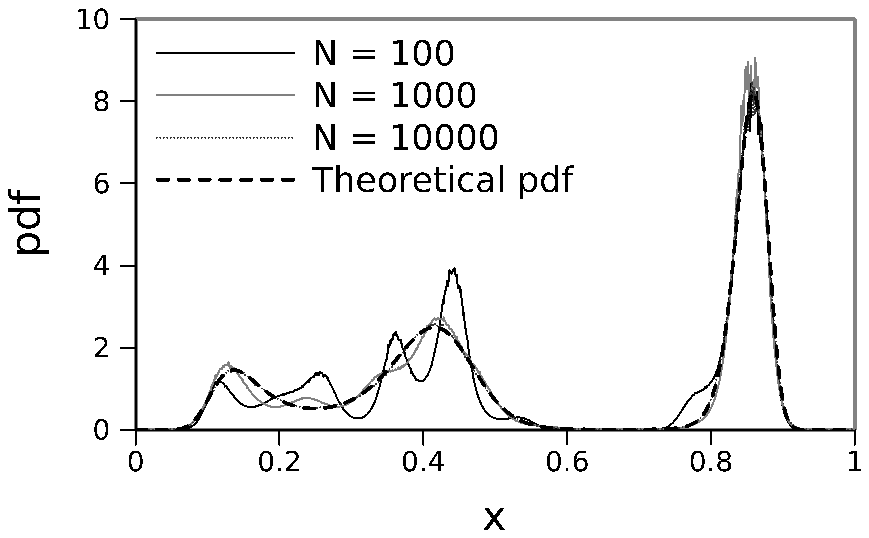}
	\includegraphics[width=.24\textwidth]{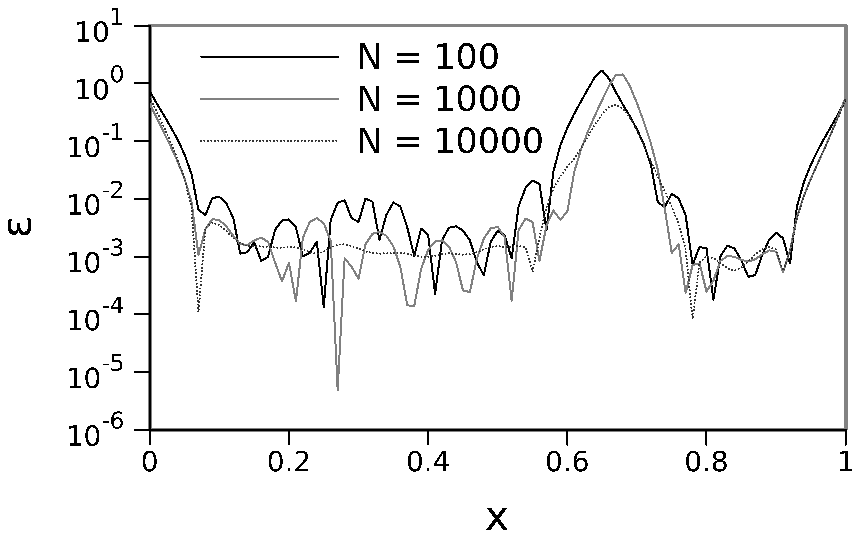}
  }
  \caption{ Pdf estimations and relative error of conditions obtained for: (first 2 graphics from left) $N=100$ and varying $\sigma$; (last 2 graphics) varying $N$ and $\sigma = \Delta X / 30$.}
\end{figure}









\end{document}